
\magnification=1200

\hsize=13.2cm
\vsize=19.5cm
\parindent=0cm   \parskip=0pt
\pageno=1

\def\ind{\hskip 1cm\relax}
\hoffset=15mm
\voffset=1cm

\ifnum\mag=\magstep1
\hoffset=0.2cm   
\voffset=-.5cm   
\fi

\pretolerance=500 \tolerance=1000  \brokenpenalty=5000

\catcode`\@=11

\font\eightrm=cmr8         \font\eighti=cmmi8
\font\eightsy=cmsy8        \font\eightbf=cmbx8
\font\eighttt=cmtt8        \font\eightit=cmti8
\font\eightsl=cmsl8        \font\sixrm=cmr6
\font\sixi=cmmi6           \font\sixsy=cmsy6
\font\sixbf=cmbx6

\skewchar\eighti='177 \skewchar\sixi='177
\skewchar\eightsy='60 \skewchar\sixsy='60

\newfam\gothfam           \newfam\bboardfam
\newfam\cyrfam

\def\tenpoint{%
  \textfont0=\tenrm \scriptfont0=\sevenrm \scriptscriptfont0=\fiverm
  \def\rm{\fam\z@\tenrm}%
  \textfont1=\teni  \scriptfont1=\seveni  \scriptscriptfont1=\fivei
  \def\oldstyle{\fam\@ne\teni}\let\old=\oldstyle
  \textfont2=\tensy \scriptfont2=\sevensy \scriptscriptfont2=\fivesy
  \textfont\gothfam=\tengoth \scriptfont\gothfam=\sevengoth
  \scriptscriptfont\gothfam=\fivegoth
  \def\goth{\fam\gothfam\tengoth}%
  \textfont\bboardfam=\tenbboard \scriptfont\bboardfam=\sevenbboard
  \scriptscriptfont\bboardfam=\sevenbboard
  \def\bb{\fam\bboardfam\tenbboard}%
 \textfont\cyrfam=\tencyr \scriptfont\cyrfam=\sevencyr
  \scriptscriptfont\cyrfam=\sixcyr
  \def\cyr{\fam\cyrfam\tencyr}%
  \textfont\itfam=\tenit
  \def\it{\fam\itfam\tenit}%
  \textfont\slfam=\tensl
  \def\sl{\fam\slfam\tensl}%
  \textfont\bffam=\tenbf \scriptfont\bffam=\sevenbf
  \scriptscriptfont\bffam=\fivebf
  \def\bf{\fam\bffam\tenbf}%
  \textfont\ttfam=\tentt
  \def\tt{\fam\ttfam\tentt}%
  \abovedisplayskip=12pt plus 3pt minus 9pt
  \belowdisplayskip=\abovedisplayskip
  \abovedisplayshortskip=0pt plus 3pt
  \belowdisplayshortskip=4pt plus 3pt
  \smallskipamount=3pt plus 1pt minus 1pt
  \medskipamount=6pt plus 2pt minus 2pt
  \bigskipamount=12pt plus 4pt minus 4pt
  \normalbaselineskip=12pt
  \setbox\strutbox=\hbox{\vrule height8.5pt depth3.5pt width0pt}%
  \let\bigf@nt=\tenrm       \let\smallf@nt=\sevenrm
  \normalbaselines\rm}

\def\eightpoint{%
  \textfont0=\eightrm \scriptfont0=\sixrm \scriptscriptfont0=\fiverm
  \def\rm{\fam\z@\eightrm}%
  \textfont1=\eighti  \scriptfont1=\sixi  \scriptscriptfont1=\fivei
  \def\oldstyle{\fam\@ne\eighti}\let\old=\oldstyle
  \textfont2=\eightsy \scriptfont2=\sixsy \scriptscriptfont2=\fivesy
  \textfont\gothfam=\eightgoth \scriptfont\gothfam=\sixgoth
  \scriptscriptfont\gothfam=\fivegoth
  \def\goth{\fam\gothfam\eightgoth}%
  \textfont\cyrfam=\eightcyr \scriptfont\cyrfam=\sixcyr
  \scriptscriptfont\cyrfam=\sixcyr
  \def\cyr{\fam\cyrfam\eightcyr}%
  \textfont\bboardfam=\eightbboard \scriptfont\bboardfam=\sevenbboard
  \scriptscriptfont\bboardfam=\sevenbboard
  \def\bb{\fam\bboardfam}%
  \textfont\itfam=\eightit
  \def\it{\fam\itfam\eightit}%
  \textfont\slfam=\eightsl
  \def\sl{\fam\slfam\eightsl}%
  \textfont\bffam=\eightbf \scriptfont\bffam=\sixbf
  \scriptscriptfont\bffam=\fivebf
  \def\bf{\fam\bffam\eightbf}%
  \textfont\ttfam=\eighttt
  \def\tt{\fam\ttfam\eighttt}%
  \abovedisplayskip=9pt plus 3pt minus 9pt
  \belowdisplayskip=\abovedisplayskip
  \abovedisplayshortskip=0pt plus 3pt
  \belowdisplayshortskip=3pt plus 3pt
  \smallskipamount=2pt plus 1pt minus 1pt
  \medskipamount=4pt plus 2pt minus 1pt
  \bigskipamount=9pt plus 3pt minus 3pt
  \normalbaselineskip=9pt
  \setbox\strutbox=\hbox{\vrule height7pt depth2pt width0pt}%
  \let\bigf@nt=\eightrm     \let\smallf@nt=\sixrm
  \normalbaselines\rm}

\def\pc#1{\bigf@nt#1\smallf@nt}         \def\pd#1 {{\pc#1} }

\def\^#1{\if#1i{\accent"5E\i}\else{\accent"5E #1}\fi}
\def\"#1{\if#1i{\accent"7F\i}\else{\accent"7F #1}\fi}

\newtoks\auteurcourant      \auteurcourant={\hfil}
\newtoks\titrecourant       \titrecourant={\hfil}

\newtoks\hautpagetitre      \hautpagetitre={\hfil}
\newtoks\baspagetitre       \baspagetitre={\hfil}

\newtoks\hautpagegauche
\hautpagegauche={\eightpoint\rlap{\folio}\hfil\the\auteurcourant\hfil}
\newtoks\hautpagedroite
\hautpagedroite={\eightpoint\hfil\the\titrecourant\hfil\llap{\folio}}

\newtoks\baspagegauche      \baspagegauche={\hfil}
\newtoks\baspagedroite      \baspagedroite={\hfil}

\newif\ifpagetitre          \pagetitretrue

\footline={\ifpagetitre\the\baspagetitre\else
\ifodd\pageno\the\baspagedroite\else\the\baspagegauche\fi\fi
\global\pagetitrefalse}

\def\raggedbottom{\topskip 10pt plus 36pt\r@ggedbottomtrue}

\def\pointir{\unskip . --- \ignorespaces}

\def\Bigbreak{\vskip-\lastskip\bigbreak}
\def\Medbreak{\vskip-\lastskip\medbreak}

\def\ctexte#1\endctexte{%
  \hbox{$\vcenter{\halign{\hfill##\hfill\crcr#1\crcr}}$}}

\long\def\ctitre#1\endctitre{%
    \ifdim\lastskip<24pt\vskip-\lastskip\bigbreak\bigbreak\fi
  		\vbox{\parindent=0pt\leftskip=0pt plus 1fill
          \rightskip=\leftskip
          \parfillskip=0pt\bf#1\par}
    \bigskip\nobreak}

\long\def\section#1\endsection{%
\vskip 0pt plus 3\normalbaselineskip
\penalty-250
\vskip 0pt plus -3\normalbaselineskip
\Bigbreak
\message{[section \string: #1]}{\bf#1\unskip}\pointir}

\long\def\sectiona#1\endsection{%
\vskip 0pt plus 3\normalbaselineskip
\penalty-250
\vskip 0pt plus -3\normalbaselineskip
\Bigbreak
\message{[sectiona \string: #1]}%
{\bf#1}\medskip\nobreak}

\long\def\subsection#1\endsubsection{%
\Medbreak
{\it#1\unskip}\pointir}

\long\def\subsectiona#1\endsubsection{%
\Medbreak
{\it#1}\par\nobreak}

\let\+=\tabalign

\def\signature#1\endsignature{\vskip 15mm minus 5mm\rightline{\vtop{#1}}}

\mathcode`A="7041 \mathcode`B="7042 \mathcode`C="7043 \mathcode`D="7044
\mathcode`E="7045 \mathcode`F="7046 \mathcode`G="7047 \mathcode`H="7048
\mathcode`I="7049 \mathcode`J="704A \mathcode`K="704B \mathcode`L="704C
\mathcode`M="704D \mathcode`N="704E \mathcode`O="704F \mathcode`P="7050
\mathcode`Q="7051 \mathcode`R="7052 \mathcode`S="7053 \mathcode`T="7054
\mathcode`U="7055 \mathcode`V="7056 \mathcode`W="7057 \mathcode`X="7058
\mathcode`Y="7059 \mathcode`Z="705A

\def\spacedmath#1{\def\packedmath##1${\bgroup\mathsurround=0pt ##1\egroup$}%
\mathsurround#1 \everymath={\packedmath}\everydisplay={\mathsurround=0pt }}

\def\nospacedmath{\mathsurround=0pt \everymath={}\everydisplay={} }

\def\decale#1{\smallbreak\hskip 28pt\llap{#1}\kern 5pt}
\def\decaledecale#1{\smallbreak\hskip 34pt\llap{#1}\kern 5pt}
\def\puce{\smallbreak\hskip 6pt{$\scriptstyle\bullet$}\kern 5pt}

\def\displaylinesno#1{\displ@y\halign{
\hbox to\displaywidth{$\@lign\hfil\displaystyle##\hfil$}&
\llap{$##$}\crcr#1\crcr}}

\def\ldisplaylinesno#1{\displ@y\halign{
\hbox to\displaywidth{$\@lign\hfil\displaystyle##\hfil$}&
\kern-\displaywidth\rlap{$##$}\tabskip\displaywidth\crcr#1\crcr}}

\def\eqalign#1{\null\,\vcenter{\openup\jot\m@th\ialign{
\strut\hfil$\displaystyle{##}$&$\displaystyle{{}##}$\hfil
&&\quad\strut\hfil$\displaystyle{##}$&$\displaystyle{{}##}$\hfil
\crcr#1\crcr}}\,}

\def\system#1{\left\{\null\,\vcenter{\openup1\jot\m@th
\ialign{\strut$##$&\hfil$##$&$##$\hfil&&
        \enskip$##$\enskip&\hfil$##$&$##$\hfil\crcr#1\crcr}}\right.}

\let\@ldmessage=\message

\def\message#1{{\def\pc{\string\pc\space}%
                \def\'{\string'}\def\`{\string`}%
                \def\^{\string^}\def\"{\string"}%
                \@ldmessage{#1}}}

\def\up#1{\raise 1ex\hbox{\smallf@nt#1}}

\def\qed{\raise -2pt\hbox{\vrule\vbox to 10pt{\hrule width 4pt
                 \vfill\hrule}\vrule}}

\def\cqfd{\unskip\penalty 500\quad\vrule height 4pt depth 0pt width
4pt\medbreak}

\def\virg{\raise .4ex\hbox{,}}   


\def\build#1_#2^#3{\mathrel{
\mathop{\kern 0pt#1}\limits_{#2}^{#3}}}

\def\boxit#1#2{%
\setbox1=\hbox{\kern#1{#2}\kern#1}%
\dimen1=\ht1 \advance\dimen1 by #1 \dimen2=\dp1 \advance\dimen2 by #1
\setbox1=\hbox{\vrule height\dimen1 depth\dimen2\box1\vrule}%
\setbox1=\vbox{\hrule\box1\hrule}%
\advance\dimen1 by .4pt \ht1=\dimen1
\advance\dimen2 by .4pt \dp1=\dimen2  \box1\relax}

\def\dateam{\the\day\ \ifcase\month\or janvier\or f\'evrier\or mars\or
avril\or mai\or juin\or juillet\or ao\^ut\or septembre\or octobre\or
novembre\or d\'ecembre\fi \ {\old \the\year}}

\def\date{\the\day\ \ifcase\month\or janvier\or f\'evrier\or mars\or
avril\or mai\or juin\or juillet\or ao\^ut\or septembre\or octobre\or
novembre\or d\'ecembre\fi \ {\old \the\year}}

\def\crog{{\vrule height 2.57mm depth 0.85mm width 0.3mm}\kern -0.36mm
[}

\def\crod{]\kern -0.4mm{\vrule height 2.57mm depth 0.85mm
width 0.3 mm}}

\def\rond{\kern 1pt{\scriptstyle\circ}\kern 1pt}

\def\hfl#1#2{\nospacedmath\smash{\mathop{\hbox to
12mm{\rightarrowfill}}\limits^{\scriptstyle#1}_{\scriptstyle#2}}}

\def\phfl#1#2{\nospacedmath\smash{\mathop{\hbox to
8mm{\rightarrowfill}}\limits^{\scriptstyle#1}_{\scriptstyle#2}}}

\def\ppav{principally polarized abelian variety}
\def\ppavs{principally polarized abelian varieties}
\def\pa{\S\kern.15em}
\def\ra{\rightarrow}
\def\lra{\longrightarrow}
\def\llra{\nospacedmath\hbox to 10mm{\rightarrowfill}}
\def\lllra{\nospacedmath\hbox to 15mm{\rightarrowfill}}
\def\saut{\vskip 5mm plus 1mm minus 2mm}

\def\Sing{\mathop{\rm Sing}\nolimits}

\def\Pic{\mathop{\rm Pic}\nolimits}

\def\codim{\mathop{\rm codim}\nolimits}

\def\cc#1{\hfill\kern .7em#1\kern .7em\hfill}

\def\dra{\ra\kern -3mm\ra}
\def\ldra{\lra\kern -3mm\ra}

\catcode`\@=12

\showboxbreadth=-1  \showboxdepth=-1

\baselineskip=14pt
\spacedmath{1.7pt}
\baspagegauche={\centerline{\tenbf\folio}}
\baspagedroite={\centerline{\tenbf\folio}}
\hautpagegauche={\hfil}
\hautpagedroite={\hfil}
\def\pa{\S\kern.15em}
\def\ra{\rightarrow}
\def\saut{\vskip 5mm plus 1mm minus 2mm}

\font\pc=cmcsc10 \rm
\parskip=1mm

\saut
\ctitre
{\pc Trisecant Lines And Jacobians, II}
\endctitre
\centerline{\pc Olivier Debarre\footnote{(*)}{\rm Partially supported by
N.S.F. Grant DMS 92-03919 and the European Science
Project ``Geometry of Algebraic Varieties", Contract no. SCI-0398-C (A).}}

\saut
{\bf 1. Introduction}

\ind Let $(X,\theta)$ be a complex \ppav. Symmetric representatives $\Theta$ of
the
polarization differ by translations by points of order $2$, hence the linear
system $|2\Theta |$ is independent on the choice of $\Theta$. It defines
a morphism $K:X\ra |2\Theta |^*$, whose image is the {\sl Kummer
variety} $K(X)$ of $X$. When $(X,\theta)$ is the Jacobian of an algebraic
curve, there are infinitely many {\it trisecants} to $K(X)$, i.e. lines
in the projective space $|2\Theta |^*$ that meet $K(X)$ in at least $3$
points. Welters conjectured in [W] that the existence of {\it one} trisecant
line to the Kummer variety should characterize  Jacobians among all
indecomposable \ppavs, thereby giving one answer to the Schottky problem.

\ind The aim of this article is to improve on the results of [D], where a
partial answer to this problem was given under additional hypotheses. More
precisely, our main theorem implies that {\it an indecomposable \ppav\
$(X,\lambda)$ is a
Jacobian if and only if there exist points $a,b,c$ of $X$ such
that:}

{\parindent=1cm\item{\rm (i)}{\it the subgroup of $X$ generated by $a-b$ and
$b-c$ is dense in $X$ ,
\item {\rm (ii)}the points $K(a)$, $K(b)$ and $K(c)$ are distinct
and collinear.\par}}

\smallskip
\ind The method is basically the same as in [D]: we prove that the existence
of one trisecant line implies the existence of a one-dimensional family
of such lines. Welters' criterion ([W]) then yields the conclusion.
\saut
{\bf 2. The set up}

\ind Let $(X,\lambda)$ be a complex indecomposable \ppav , let
$\Theta$ be a symmetric representative of the polarization and let $K:X\ra
|2\Theta |^*$ be the Kummer morphism. Let $\theta$ be a non-zero section of
${\cal O}_X
(\Theta)$. For any $x\in X$, we write $\Theta_x$ for the divisor $\Theta +x$
and $\theta _x$
for the section $z\mapsto \theta (z-x)$ of ${\cal O}_X (\Theta _x)$. If $a,b$
and $c$ are points of $X$, it is classical that the points $K(a), K(b)$
and $K(c)$ are collinear if and only if there exist complex numbers $\alpha ,
\beta$ and $\gamma$ not all zero such that:
$$
\alpha \theta _a\theta _{-a} + \beta \theta _b\theta _{-b} + \gamma \theta
_c\theta _{-c} = 0\ .
$$
Following Welters, we consider the set:
$$
V_{a,b,c}=2\ \{ \zeta \in X \mid K(\zeta +a),K(\zeta +b),K(\zeta +c)\
\hbox {are collinear}\}
$$
endowed with its natural scheme structure. By [W], theorem 0.5,  $(X,\lambda)$
{\it is a
Jacobian if and only if there exist points} $a,b$ {\it and} $c$ {\it such that}
$\dim
V_{a,b,c}>0$. This condition is equivalent to the existence of a sequence $\{
D_n\}_{n>0}$ of
constant vector fields on $X$ and of a formal curve $\zeta (\epsilon ) = \zeta
(0)+{1\over
2}D(\epsilon )$ with $D(\epsilon )=\sum_{n>0} D_n \epsilon ^n$, contained in
$V_{a,b,c}$. This
in turn is equivalent to a relation of the type:
$$
\alpha (\epsilon ) \theta _{a+\zeta (\epsilon)}\theta _{-a-\zeta (\epsilon)} +
\beta (\epsilon ) \theta _{b+\zeta (\epsilon)}\theta _{-b-\zeta (\epsilon)} +
\gamma (\epsilon ) \theta _{c+\zeta (\epsilon)}\theta _{-c-\zeta (\epsilon)} =
0\ ,
$$
where $\alpha (\epsilon ),
\beta (\epsilon )$ and $\gamma (\epsilon )$ are relatively prime elements of
${\bf C} [[\epsilon ]]$.
\saut
{\bf 3. The case of a degenerate trisecant}

\ind In this section, we prove the following result:
\medskip

{\pc Theorem} 3.1.-- {\it Let $(X,\lambda)$ be a complex indecomposable \ppav ,
let $\Theta$ be
a symmetric representative of the polarization and let $K:X\ra
|2\Theta |^*$ be the Kummer morphism. Assume that
there exist two points $u$ and $v$ of $X$ such that:

{\parindent=1cm\item{\rm (i)}{\it the points $K(u)$ and $K(v)$ are distinct and
non-singular on
$K(X)$ and the line that joins them is tangent to $K(X)$ at $K(u)$,
\item {\rm (ii)}$\displaystyle{\rm codim} _X \bigcap_{s \in {\bf Z}}
\Theta_{2su} > 2\ $.}\par}
Then $(X,\lambda)$ is isomorphic to the Jacobian of a smooth algebraic curve.}
\medskip
\ind Note that condition (ii) in the theorem holds when $u$ generates $X$.
\smallskip
{\bf Proof.} As explained in section 2, it is enough to prove that the scheme
$V_{u,-u,v}$
has positive dimension at $0$: we look for a sequence $\{ D_n\}_{n>0}$ of
constant vector fields
on $X$ with $D_1\not = 0$ and relatively prime elements $\alpha (\epsilon ),
\beta (\epsilon )$ and $\gamma (\epsilon )$ of
${\bf C} [[\epsilon ]]$ such that:
$$
\alpha (\epsilon ) \theta _{u+{1\over 2} D(\epsilon)}\theta _{-u-{1\over 2}
D(\epsilon)} +
\beta (\epsilon ) \theta _{-u+{1\over 2} D(\epsilon)}\theta _{u-{1\over 2}
D(\epsilon)} +
\gamma (\epsilon ) \theta _{v+{1\over 2} D(\epsilon)}\theta _{-v-{1\over 2}
D(\epsilon)} =
0\ ,\leqno (3.2)
$$
with $D(\epsilon )=\sum_{n>0} D_n \epsilon ^n$. This is nothing but equation
(1.4) from [D].
It follows from {\it loc.cit.} that we may assume:
$$
\alpha (\epsilon) = 1+\sum _{n>0} \alpha _n
\epsilon ^n \quad ,\quad \beta (\epsilon) =-1\quad ,\quad \gamma (\epsilon)
=\epsilon\quad .
$$
Write the left-hand-side of (3.2) as $\sum _{n\geq 0} P _n \epsilon ^n$, where,
for any
$n\geq 0$, $P_n$ is a section of ${\cal O}_X(2\Theta)$. One has $P_0=0$ and:
 $$
P_1= \alpha_1\theta _u \theta _{-u} + \theta _u D_1\theta _{-u} -\theta _{-u}
D_1\theta _u +
 \theta _v\theta _{-v}\ . \leqno (3.3)
$$
\ind As explained in {\it loc.cit.}, hypothesis (i) in the theorem is
equivalent to the
vanishing of $P_1$ for a suitable $D_1$ tangent at $K(u)$ to the line that
joins $K(u)$ and
$K(v)$, and a suitable $\alpha_1$. In general, note that $P_n$ depends only on
$\alpha
_1,\ldots,\alpha_n$ and $D_1,\ldots,D_n$. Knowing that $P_1$ vanishes, we need
to construct
a sequence $\{ D_n\}_{n>0}$ of constant vector fields on
$X$ and  a sequence  $\{ \alpha_n\}_{n>1}$ of complex numbers such that $P_n$
vanishes for all
positive integers $n$.

\ind We proceed by induction: let $n$ be an integer $\geq 2$ and assume that
$\alpha
_1,\ldots,\alpha_{n-1}$ and $D_1,\ldots,D_{n-1}$ have been constructed so that
$P_1=\cdots=P_{n-1}=0$. We want to find a complex number $\alpha_n$ and a
tangent vector $D_n$
such that $P_n$ vanishes on $X$. By lemma 1.8 from {\it loc.cit.}, it is enough
to show that
the restriction of $P_n$ to the scheme $\Theta_u \cap \Theta_{-u}$ (which
depends only on
$\alpha _1,\ldots,\alpha_{n-1}$ and $D_1,\ldots,D_{n-1}$) vanishes. It is
convenient to set:
$$
R(z,\epsilon)=P(z+{1\over 2}D(\epsilon),\epsilon) = \sum _{n>0}
R_n(z)\epsilon^n\ . $$
\ind Our induction hypothesis can be rewritten as $R_1=\cdots=R_{n-1}=0$ and
$P_n=R_n$. Therefore, we need to prove that $R_n$ vanishes on the scheme
$\Theta_u \cap
\Theta_{-u}$. We begin by proving a few identities independent on the induction
hypothesis.
Note that: $$
R(z,\epsilon)=
\alpha (\epsilon ) \theta (z-u)\theta(z+u+D(\epsilon)) -
\theta (z+u)\theta(z-u+D(\epsilon)) + \epsilon  \theta
(z-v) \theta(z+v+D(\epsilon)) \ .\leqno (3.4)
$$
It follows that for any $z$ in $\Theta_u$, one has:
$$\leqalignno{
\qquad\  & R(z,\epsilon)
=-\theta(z+u)\theta(z-u+D(\epsilon ))+\epsilon\theta(z-v)\theta(z+v+D(\epsilon
)) & (3.5) \cr
\qquad\  & R(z-2u,\epsilon)
=\alpha(\epsilon)\theta(z-3u)\theta(z-u+D(\epsilon
))+\epsilon\theta(z-2u-v)\theta(z-2u+v+D(\epsilon )) & (3.6) \cr
\qquad\  & R(z-u+v,\epsilon)
=\alpha(\epsilon)\theta(z-2u+v)\theta(z+v+D(\epsilon
))-\theta(z+v)\theta(z-2u+v+D(\epsilon
))\ \ . & (3.7) \cr }
$$
Moreover, since $P_1$ and its translate by $2u$ both vanish, formula (3.3)
yields, for any
$z$ in $\Theta_u$:
$$\leqalignno{
\theta(z+u)D_1\theta(z-u) - \theta(z+v)\theta(z-v) & =0 & (3.8) \cr
\theta(z-3u)D_1\theta(z-u) + \theta(z-2u+v)\theta(z-2u-v) & =0. & (3.9) \cr
}$$
\saut
The following result is the main technical step of the proof.
\medskip
{\pc Lemma} 3.10.-- {\it For any point $z$ in $\Theta_u$, one has:}
\nospacedmath
$$\displaylines{
\qquad \alpha(\epsilon)R(z,\epsilon)\theta(z-3u)D_1\theta(z-u)
+ R(z-2u,\epsilon)\theta(z+u)D_1\theta(z-u)  \hfill\cr
\hfill {}+\epsilon
R(z-u+v,\epsilon)\theta(z-v)\theta(z-2u-v)=0\ . \qquad\cr}
$$
\spacedmath{1.7pt}
\vskip -.3truecm
{\bf Proof.} By (3.5), (3.6) and (3.7), the left-hand-side of the expression in
the lemma is
equal to:
\nospacedmath
$$\eqalign{
& -\alpha(\epsilon)\ \theta(z-u+D(\epsilon))\ \theta(z+u)\ \theta(z-3u)\
D_1\theta(z-u) \cr
  {} & +\epsilon\alpha(\epsilon)\ \theta(z+v+D(\epsilon))\ \theta(z-v)\
\theta(z-3u)\
D_1\theta(z-u) \cr {} & +\alpha(\epsilon)\ \theta(z-u+D(\epsilon))\
\theta(z-3u)\ \theta(z+u)\
D_1\theta(z-u) \cr {} & +\epsilon\ \theta(z-2u+v+D(\epsilon))\ \theta(z-2u-v)\
\theta(z+u)\
D_1\theta(z-u) \cr {} & +\epsilon\alpha(\epsilon)\ \theta(z+v+D(\epsilon))\
\theta(z-2u+v)\
\theta(z-v)\ \theta(z-2u-v) \cr
{} & -\epsilon\ \theta(z-2u+v+D(\epsilon))\ \theta(z+v)\ \theta(z-v)\
\theta(z-2u-v)\ .\cr }
$$
\spacedmath{1.7pt}
In this sum, the first and third terms cancel out; the second and fifth cancel
out by
(3.9) and the fourth and sixth by (3.8). Hence the sum vanishes.\cqfd
\saut
\ind Recall that we are assuming $R_1=\cdots=R_{n-1}=0$. Since
$\alpha(\epsilon)\equiv 1$
modulo $\epsilon$, the identity of the lemma taken modulo $\epsilon^{n+1}$
yields:
$$
 D_1\theta(z-u) \bigl\lbrack
R_n(z)\theta(z-3u)+R_n(z-2u)\theta(z+u)\bigr\rbrack = 0\ .\leqno \qquad \forall
z\in \Theta_u
 $$
Since $\Theta_u$ is integral and $D_1$ is non-zero, we get:
$$
R_n(z)\theta(z-3u)+R_n(z-2u)\theta(z+u) = 0\ .\leqno (3.11)\qquad\forall z\in
\Theta_u
 $$
\medskip
{\pc Lemma} 3.12.-- {\it If $F$ is a section of a line bundle on $X$ such that
$R_nF$
vanishes on the scheme $\Theta_u \cap\Theta_{-u}$, then, for any integer $s$,
the section $R_n
F_{2su}$ also vanishes on the scheme $\Theta_u  \cap\Theta_{-u}$.}
\smallskip
{\bf Proof.} It is enough to prove that $R_n
F_{2u}$ vanishes on $\Theta_u \cap \Theta_{-u}$ since $R_n
F_{-2u}$ is its image by the involution $x\mapsto -x$. There exist two sections
$A$ and $B$
such that $R_n F=A\theta_u+B\theta_{-u}$. It follows that $(R_n)_{2u}
F_{2u}\equiv
A_{2u}\theta_{3u}\pmod {\theta_u}$. Multiplying (3.11) by $F_{2u}$, we get:
$$
R_n\ F_{2u}\ \theta_{3u}+A_{2u}\ \theta_{3u}\ \theta_{-u}\equiv 0\eqno
\pmod {\theta_u}\ .\qquad
 $$
Since $2u\not = 0$ (because $K(u)$ is non-singular) and $\Theta_u$ is integral,
$\theta_{3u}$ is
not a zero divisor modulo $\theta_u$ and we get
$
R_n\ F_{2u}\equiv 0\ ({\rm mod}\ (\theta_u,\theta_{-u}))
 $.\cqfd
\smallskip
\ind  The lemma immediately yields that  $R_n\theta_{u+2su}$ vanishes on
$\Theta_u
\cap\Theta_{-u}$ for all $s$. Hypothesis (ii) in the theorem guaranties that
the scheme
$\bigcap_{s \in {\bf Z}} \Theta_{u+2su}$ has codimension $>2$ in $X$. Since the
scheme $\Theta_u
\cap \Theta_{-u}$ has pure codimension $2$, it follows that
 $R_n$ vanishes there. This concludes the proof of the theorem.\cqfd
\smallskip
\ind Theorem 6.2 below shows that the conclusion of  theorem 3.1 still holds
with a
hypothesis slightly different from (ii).
 \saut
 {\bf 4. The case of a non-degenerate trisecant}

\ind In this section, we prove, under an extra hypothesis, that the existence
of a
non-degenerate trisecant line implies the existence of a degenerate trisecant
of the type
studied in section 3. We prove:
\medskip

{\pc Theorem} 4.1.-- {\it Let $(X,\lambda)$ be an indecomposable \ppav , let
$\Theta$ be a
symmetric representative of the polarization and let $K:X\ra
|2\Theta |^*$ be the Kummer morphism. Assume that
there exist points $a,b$ and $c$ of $X$ such that:

{\parindent=1cm\item{\rm (i)}{\it the points $K(a), K(b)$ and $K(c)$ are
distinct and collinear,
\item {\rm (ii)}${\rm codim}_X \displaystyle\bigcap_{\scriptstyle p,q,r \in
{\bf
Z}\atop\scriptstyle p+q+r=0} \Theta_{pa+qb+rc} > 2\ $.}\par}
\smallskip
\ind Then $(X,\lambda)$ is isomorphic to the Jacobian of a
smooth algebraic curve.}
\medskip
\ind Note that condition (ii) in the theorem holds when $a-b$ and $b-c$
together generate $X$.
\smallskip
{\bf Proof.} Instead of proving that
$V_{a,b,c}$ has positive dimension at $0$, we
will proceed as follows. As explained in section 2, condition (i) translates
into the existence
of non-zero complex numbers $\alpha, \beta$ and $\gamma$ such that
$$
\alpha \theta _a\theta _{-a} + \beta \theta _b\theta _{-b} + \gamma \theta
_c\theta _{-c} = 0\ .
\leqno (4.2)
$$
For any   $x$ in $X$, we will write $P^x$ for $\theta_{a+b+c}\theta_{-x}$. Our
first
aim is to show that $P^c$ vanishes on the scheme $\Theta_a\cap\Theta_b$.
\medskip
{\pc Lemma} 4.3.-- {\it One has:}
$$
P^c_{a-b}\theta^{}_b+P^c\theta_{2a-b}\equiv 0
\eqno ({\rm mod}\ \theta_a)\ .\qquad
$$
\smallskip
{\bf Proof.} Equation (4.2) and its translates by $(a+c)$ and $(a-b)$ yield,
modulo $\theta_a$:
$$\eqalign{
\beta\theta _b\theta _{-b} + \gamma \theta
_c\theta _{-c} & \equiv 0 \cr
\alpha \theta _{2a+c}\theta _c + \beta \theta _{a+b+c}\theta _{a-b+c}  & \equiv
0 \cr
\alpha \theta _{2a-b}\theta _{-b} + \gamma \theta _{a-b+c}\theta _{a-b-c}  &
\equiv 0\ . \cr
}$$
It follows that, still modulo $\theta_a$:
$$\eqalign{
\alpha\beta\theta _{-b}\ (  P^c_{a-b}\theta^{}_b & +P^c\theta_{2a-b}) \cr
& \equiv \theta _{2a+c}\theta _{a-b-c}\ (-\alpha\gamma\theta
_c\theta _{-c})+\theta _{a+b+c}\theta
_{-c}\ (-\beta\gamma\theta
_{a-b+c}\theta _{a-b-c}) \cr
 & \equiv -\gamma\theta
_{a-b-c}\theta_{-c}(\alpha\theta_{2a+c}\theta_c+\beta\theta_{a+b+c}\theta_{a-b+c}) \cr
& \equiv\ 0\ . \cr
}$$
Since $\Theta$ is integral and $-b\not = a$, the section $\theta_{-b}$ is not a
zero
divisor modulo $\theta_a$ and the lemma follows.\cqfd

 \medskip
{\pc Lemma} 4.4.-- {\it If $F$ is a section of a line bundle on $X$ such that
$P^cF$
vanishes on the scheme $\Theta_a\cap \Theta_b$, then, for any integer $s$, the
section $P^c
F_{s(a-b)}$ also vanishes on the scheme $\Theta_a\cap \Theta_b$.}
\smallskip
{\bf Proof.} Since $a$ and $b$ play the same role, it is enough to prove that
$P^c
F_{a-b}$ vanishes on $\Theta_a\cap \Theta_b$. Let $A$ and $B$ be two sections
such that
$P^c F=A\theta_a+B\theta_b$. Then $P^c_{a-b} F^{}_{a-b}\equiv
A_{a-b}\theta_{2a-b}\ ({\rm mod}\ \theta_a)$. Using lemma 4.3, we get:
$$
P^c\ F_{a-b}\ \theta_{2a-b}+A_{a-b}\ \theta_{2a-b}\ \theta_{b}\equiv 0\eqno
({\rm mod}\ \theta_a)\ .\qquad
$$
Since $(2a-b)\not = a$ and $\Theta_a$ is irreducible, we can divide out by
$\theta_{2a-b}$,
and the lemma is proved.\cqfd
\medskip
{\pc Lemma} 4.5.-- {\it If $F$ is a section of a line bundle on $X$, then
$P^cF$
vanishes on the scheme $\Theta_a \cap\Theta_b$ if and only if $P^bF$
vanishes on the scheme $\Theta_a\cap \Theta_c$.}
\smallskip
{\bf Proof.} Write:
$$
\theta_{a+b+c}\ \theta_{-c}\ F \equiv A\ \theta_b
\eqno ({\rm mod}\ \theta_a)\ . \qquad
$$
Then,  on $\Theta_a$:
$$\eqalign{
\gamma\ A\ \theta_b\ \theta_c & \equiv \theta _{a+b+c}\ \gamma\ \theta
_c\ \theta _{-c}\ F \cr
& \equiv -\ \theta _{a+b+c}\ \beta\ \theta
_b\ \theta _{-b}\ F \cr
& = -\ \beta \ P^b\ F\ \theta_b\ , \cr
}$$
where we used (4.2). Since $\Theta_a$ is irreducible and $a\not = b$, the lemma
is
proved.\cqfd
\medskip
\ind We now combine the last two lemmas to get, for all integers $r$ and $s$:
\medskip
\vbox{
\halign{\quad\qquad#&\quad\quad\quad#&\quad\quad#&\quad\quad#\cr
 & $\quad\qquad\qquad P^c\ F$ & $\equiv\ 0$ & $\quad \bigl( {\rm mod}\
(\theta_a,\theta_b)\bigr) $ \cr $\Longrightarrow$ & $\quad\qquad\qquad P^b\ F$
& $\equiv\ 0$ &
$\quad \bigl( {\rm mod}\ (\theta_a,\theta_c)\bigr) $ \cr $\Longrightarrow$ & $\
\ \qquad P^b\
F_{r(a-c)}$ & $\equiv\ 0$ & $\quad \bigl( {\rm mod}\ (\theta_a,\theta_c)\bigr)
$ \cr
$\Longrightarrow$ & $\ \ \qquad P^c\ F_{r(a-c)}$ & $\equiv\ 0$ & $\quad \bigl(
{\rm mod}\
(\theta_a,\theta_b)\bigr)$ \cr $\Longrightarrow$ & $P^c\ F_{r(a-c)+s(a-b)}$ &
$\equiv\ 0$ &
$\quad \bigl( {\rm mod}\ (\theta_a,\theta_b)\bigr)\ .$ \cr }}

\medskip
It follows in particular that $P^c\theta_{a+r(a-c)+s(a-b)}$ vanishes on
$\Theta_a\cap\Theta_b$ for all integers $r$ and $s$. Hypothesis (ii) in the
theorem then
implies:
$$P^c\ {\rm vanishes\ on\ the\ scheme}\ \Theta_a\cap\Theta_b\leqno (4.6)$$
(hence also $P^a$ on $\Theta_b\cap\Theta_c$
and $P^b$  on $\Theta_c\cap\Theta_a$).

\ind Let $u$ be any point of $X$ such that $2u=a-b$ and set $v=u-a-c$.
Translating (4.6)
by $(-u-b)$, we get that $\theta_v\theta_{-v}$ vanishes on
$\Theta_u\cap\Theta_{-u}$. As
explained in [D], this is equivalent to the existence of a complex number
$\alpha_1$ and a
tangent vector $D_1$ to $X$ such that:
 $$
\alpha_1\theta _u \theta _{-u} + \theta _u D_1\theta _{-u} -\theta _{-u}
D_1\theta _u +
 \theta _v\theta _{-v}=0\ . \leqno (4.7)
$$
In other words, the line that joins $K(u)$ and $K(v)$ is tangent to $K(X)$ at
$K(u)$. Note
that we cannot apply theorem 3.1 directly, since hypothesis (ii) is not
satisfied. However,
we will still follow the same method, i.e{.}{\ }we will show that the scheme
$V_{a,b,-c}$ (which
is a translate of $V_{u,-u,v}$) has positive dimension at $(-a-b)$, but we will
need to prove at
the same time that $V_{a,-b,c}$ has positive dimension at the point $(-a-c)$.

\ind Let $n$ be an integer $\geq 1$. As in section 3, the scheme $V_{a,b,-c}$
contains a scheme isomorphic to ${\bf C}[\epsilon]/\epsilon^{n+1}$ and
concentrated at $(-a-b)$
if and only if one can find complex numbers $\alpha_1,\ldots,\alpha_n$ and
tangent vectors
$D_1,\ldots,D_n$ such that $R_1,\ldots ,R_n$, defined in section 3, vanish
($\alpha_1$ and $D_1$
are the same as in (4.7), and $R_1$ is the left-hand-side of that equation).
Similarly, the
scheme $V_{a,-b,c}$ contains a scheme isomorphic to ${\bf
C}[\epsilon]/\epsilon^{n+1}$ and
concentrated at $(-a-c)$ if and only if there exist complex numbers
$\alpha'_1,\ldots,\alpha'_n$ and tangent vectors $D'_1,\ldots,D'_n$ such that
$R'_1,\ldots
,R'_n$ vanish.

\ind We proceed as in section 3: let $n$ be an integer $\geq 2$ and assume that
$\alpha
_1,\ldots,\alpha_{n-1},$ $\alpha'
_1,\ldots,\alpha'_{n-1}$ and $D_1,\ldots,D_{n-1},D'_1,\ldots,D'_{n-1}$ have
been constructed so
that $R_1,\ldots,R_{n-1},$ $R'_1,\ldots,R'_{n-1}$ vanish on $X$. As in the
proof of
theorem 3.1, it is enough to show that the restriction of $R_n$ to the scheme
$\Theta_a\cap\Theta_b$ (which depends only on $\alpha_1,\ldots,\alpha_{n-1}$
and
$D_1,\ldots,D_{n-1}$), and the restriction of $R'_n$ to
$\Theta_a\cap\Theta_c$ (which depends only on $\alpha'_1,\ldots,\alpha'_{n-1}$
and
$D'_1,\ldots,D'_{n-1}$) both vanish. By lemma 3.12, we have:
\vskip -5mm plus .5mm minus .5mm
\nospacedmath $$\normalbaselineskip=.6truecm\matrix{
R_n\theta_{a+s(a-b)}
{\rm\ \ vanishes\ on\ the\ scheme\ \ } \Theta_a\cap\Theta_b \cr
R_n\theta_{a+s(a-c)}
{\rm\ \ vanishes\ on\ the\ scheme\ \ } \Theta_a\cap\Theta_c \cr
}\leqno (4.8)$$
\spacedmath{1.7pt}
 for all integers $s$.
\medskip
{\pc Lemma} 4.9.-- {\it
For any integer $s$ such that $0<s<n$, one has $\beta^sD'_s=(-\gamma)^sD_s$.}
\smallskip
{\bf Proof.} Formula (3.4) translates into:
\nospacedmath
$$\displaylines{
(4.10)\qquad R(z,\epsilon)  =
\alpha (\epsilon ) \theta (z-a)\theta(z-b+D(\epsilon)) -
\theta (z-b)\theta(z-a+D(\epsilon))  \hfill\cr
\hfill {}+ \epsilon  \theta
(z+c) \theta(z-a-b-c+D(\epsilon))\ \,\,  \quad\cr
(4.11)\quad R'(z,\epsilon)  =
\alpha' (\epsilon ) \theta (z-a)\theta(z-c+D'(\epsilon)) -
\theta (z-c)\theta(z-a+D'(\epsilon))  \hfill\cr
\hfill {}+ \epsilon  \theta
(z+b) \theta(z-a-b-c+D'(\epsilon))\ .   \quad\cr
}
$$
\spacedmath{1.7pt}
Using (4.2), we get:
$$
\gamma R_1\theta_c + \beta R'_1\theta_b \equiv -\theta_b\theta_c(\gamma
D_1+\beta
D'_1)\theta_a \eqno ({\rm mod}\ \theta_a)\ .\qquad
$$
Since $R_1$ and $R'_1$ vanish, we  get $\gamma D_1+\beta
D'_1=0$. We complete the proof by induction on $s$. Assume that $s<n$ and that
$\beta^tD'_t=(-\gamma)^tD_t$ whenever $0<t<s$. This is equivalent to
$D'(\beta\epsilon)\equiv D(-\gamma\epsilon) $ (mod $\epsilon^s$). Using again
(4.10) and
(4.11), we get:
$$
\qquad R(\cdot,-\gamma\epsilon)\theta_c -R'(\cdot,\beta\epsilon)\theta_b \equiv
-\theta_b\theta_c\bigl( (-\gamma)^s D_s-\beta^s D'_s\bigr) \theta_a
\qquad\bigl( {\rm
mod}\ (\theta_a,\epsilon^{s+1})\bigr)\ .\qquad \leqno (4.12)
$$
Since $R$ and $R'$ vanish modulo $\epsilon^{s+1}$, we get $(-\gamma)^s
D_s-\beta^s D'_s=0$ and the lemma follows.\cqfd
\medskip
{\pc Lemma} 4.13.-- {\it Let $F$ be a section of a line bundle on $X$. Then
$R_nF$ vanishes on
the scheme $\Theta_a\cap\Theta_b$ if and only if $R'_nF$ vanishes on the scheme
$\Theta_a\cap\Theta_c$.}  \smallskip  {\bf Proof.} Our assumption is that
$R_nF\equiv A\theta_b$
(mod $\theta_a$). Formula (4.12) is still valid for $s=n$ and reads
$(-\gamma)^nR_n\theta_c
-\beta^nR'_n\theta_b \equiv -\theta_b\theta_c\bigl( (-\gamma)^n D_n-\beta^n
D'_n\bigr)
\theta_a$ (mod $\theta_a$). Multiplying this congruence by $F$, we get:
$$(-\gamma)^nA\theta_b\theta_c -\beta^nR'_nF\theta_b \equiv
-\theta_b\theta_cF\bigl(
(-\gamma)^n D_n-\beta^n D'_n\bigr) \theta_a
\eqno
\pmod {\theta_a}\ .\qquad $$
Since $\Theta_a$ is irreducible
and  $a\ne b$, one can divide out by $\theta_b$. This finishes the proof of the
lemma.\cqfd \medskip
\ind
The argument used right after the proof of lemma 4.5 and (4.8) immediately
yield that if $R_nF$
vanishes on $\Theta_a\cap\Theta_b$, then $R_n\theta_{a+r(a-c)+s(a-b)}$ vanishes
on
$\Theta_a\cap\Theta_b$ for all integers $r$ and $s$. Hypothesis (ii) then
implies that $R_n$
vanishes on $\Theta_a\cap\Theta_b$, which concludes the proof of the
theorem.\cqfd
 \saut
{\bf 5. The most degenerate case: the K-P equation}

\ind On the Jacobian of a smooth projective curve, there exists a
$2$-dimensional family of
trisecant lines. If one lets the points of contact tend to a singular point of
the Kummer
variety in a suitable way, one gets that the theta function $\theta$ of a
Jacobian satisfies
the Kadomcev-Petviashvili (or K-P) equation:
\nospacedmath
$$\displaylines{
\qquad D_1^4\theta\cdot\theta-4D_1^3\theta\cdot
D_1\theta+3(D_1^2\theta)^2-3(D_2\theta)^2
\hfill \cr
\hfill
{}+3D^2_2\theta\cdot\theta+3D_1\theta\cdot
D_3\theta-3D_1D_3\theta\cdot\theta+d_3\theta^2\ =\
0\ ,\qquad \cr }$$
\spacedmath{1.7pt}for some constant vector fields $D_1,D_2,D_3$ and complex
number $d_3$.
Shiota proved in [S] that Jacobians are characterized among indecomposable
\ppavs\  by this
property. His proof was later simplified by Arbarello and De Concini in [AC],
but they still
had to rely on the technical but crucial lemma 7 from [S]. In this section, we
will try to
relate our approach to the (discrete) trisecant conjecture to the techniques
used in those two
articles.

\ind Let $(X,\lambda)$ be an indecomposable \ppav\  which satisfies the K-P
equation, and let
$\Theta$ be a symmetric representative of the polarization. The article [AC]
reduces the proof
that $(X,\lambda)$ is a Jacobian to showing that each element of a certain
sequence
$\{P_n\}_{n>2}$ of sections of ${\cal O}_X(2\Theta)$ vanishes, where $P_3=0$ is
the K-P
equation, and where each $P_n$ depends on tangent vectors $D_1,\ldots,D_n$ and
complex numbers
$d_3,\ldots,d_n$ (where $D_1,D_2,D_3$ and $d_3$ are the same as in the K-P
equation). Let $n$
be an integer $>3$. As we did in sections 3 and 4, we assume that
$D_1,\ldots,D_{n-1}$ and
$d_3,\ldots,d_{n-1}$ are such that $P_3,\ldots,P_{n-1}$ vanish; it  was shown
in {\it loc.cit.} that it is enough to prove that the restriction of $P_n$ to
the scheme
$D_1\Theta$ defined by $\theta=D_1\theta=0$ vanishes. This restriction depends
only on
$D_1,\ldots,D_{n-1}$ and  $d_3,\ldots,d_{n-1}$.  \medskip
{\pc Step} 1.-- {\it The following
relation holds on $\Theta$: $$ D_1P_n\cdot D_1\theta -P_n\cdot D_1^2\theta=0\
.$$}
\vskip -.4truecm
{\bf Proof.} This identity should be thought of as the analog of (3.11) and
ought to
follow from a relation similar to 3.10, but I was unfortunately unable to find
it. It is
however a consequence of the following equality, valid on $X$ and proved
algebraically in [AC]
under the assumption that $P_3,\ldots,P_{n-1}$ vanish:
$$
\theta^4\Bigl( {1\over3}D_1^4+D_2^2-D_1D_3+4D_1\bigl( (D_1^2{\rm
log}\theta)\cdot
D_1\bigr)\Bigr)P_n=0\ . $$
Indeed, a direct calculation shows that the left-hand-side is equal to:
$$
8\ {(D_1\theta)^2\over \theta}\ (P_n\cdot D_1^2\theta - D_1P_n\cdot D_1\theta)\
+\
 {\rm regular\ fonction\ }.
$$
This finishes the proof.\cqfd
\medskip
\ind Another interpretation of this first step is that it is equivalent to the
vanishing of
the residue necessary to solve locally, outside of $D_1\Theta$, the equation
$\displaystyle D_1h={P_n\over \theta^2}$. This is the method used in [AC].
\medskip
{\pc Step} 2.-- {\it Let $F$ be a section of a line bundle on $X$. If $P_nF$
vanishes on
$D_1\Theta$, then $P_nD_1F$ vanishes on
$D_1\Theta$.}  \smallskip  {\bf Proof.} Write $P_nF=A\theta+BD_1\theta$ and
take the
$D_1$--derivative: $$
D_1P_n\cdot F+P_n\cdot D_1F=D_1A\cdot \theta+A\cdot D_1\theta+D_1B\cdot
D_1\theta+B\cdot
D_1^2\theta\ .  $$
Multiply by $P_n$ and use step 1 and the equality $P_nF=A\theta+BD_1\theta$ to
get:
$$
D_1P_n\cdot B\cdot D_1\theta+P_n^2\cdot D_1F\equiv(A+D_1B)\cdot P_n\cdot
D_1\theta+B\cdot
D_1P_n\cdot D_1\theta\eqno ({\rm mod}\ \theta)\ ,\qquad $$
hence
$$
P_n\bigl(P_n\cdot D_1F-(A+D_1B)\cdot D_1\theta\bigr)\equiv 0\eqno ({\rm mod}\
\theta)\ .\qquad
$$
Since $\Theta$ is irreducible, one of the two factors vanishes on $\Theta$,
which proves
our contention.\cqfd
\medskip
\ind It follows that for any positive integer $s$, the section $P_n\cdot
D_1^s\theta$ vanishes on $D_1\theta$. If the codimension of $\Sigma
=\bigcap_{s\geq
0}D_1^s\Theta$ in $X$ is $>2$, the section $P_n$ vanishes on $D_1\Theta$. This
yields a
completely algebraic proof of the Novikov conjecture, under the additional
hypothesis
${\rm codim}_X\bigcap_{s\geq 0}D_1^s\Theta > 2$, which slightly improves on the
main theorem of [AC]. See also the end of the next section for another
algebraic proof with a
slightly different hypothesis.

\ind To get rid of this extra assumption, one needs lemma 7
from [S], which shows essentially that the K-P equation implies that  $\Sigma$
is
equal to  $\bigcap_{r,s,t\geq 0} D_1^rD_2^sD_3^t\Theta$ (as schemes, although
we will only use the
equality as sets). Following Shiota, we work on a desingularization $\widetilde
X$ of the
blow-up of the latter scheme in $X$: the K-P equation still makes sense on
$\widetilde X$, because
the vector fields $D_1,D_2$ and $D_3$ are still defined there; one checks that
the strict transform
$\widetilde{\Theta}$ of $\Theta$ is irreducible, and lemma 7 of [S] applies on
$\widetilde X$
to show that  $\bigcap_{s\geq 0}D_1^s\widetilde{\Theta}$ is equal to
$\bigcap_{r,s,t\geq 0} D_1^rD_2^sD_3^t\widetilde{\Theta}$, hence is empty. The
above ideas
applied on $\widetilde X$ then give a proof of the Novikov conjecture.

\ind Note   that  $\bigcap_{s\geq
0}D_1^s\Theta$ is  a posteriori actually empty, whereas in the other cases, it
may
happen that   $\bigcap_{s \in {\bf Z}} \Theta_{2su}$ or $\bigcap_{p,q \in {\bf
Z}}
\Theta_{p(a-c)+q(b-c)}$ are not empty (if for example $u$, or $a-c$ and $b-c$,
are torsion).

\saut
{\bf 6. Complements}

\ind In this short section, we will indicate how to combine the techniques used
here whith those of
[D] to get results in the degenerate cases when the theta divisor in not too
singular. We will use
the following lemma, inspired by proposition 2.6 in [D].
\medskip

{\pc Lemma} 6.1.-- {\it Let $(X,\lambda)$ be an indecomposable \ppav and let
$\Theta$ be a
 representative of the polarization.  Let $x$ be a
non-torsion element of $X$ and assume that $Z$ is a component of
$\Theta\cap\Theta_x$
such that $Z_{\rm red}$ is contained in $\bigcap_{s \in {\bf Z}} \Theta_{sx}$.
Assume that $\codim
_X (Z\cap\Sing\Theta)>3$.Then $Z$ is reduced.} \smallskip

{\bf Proof.}
Since $\bigl( Z_{\rm red}\bigr)_{sx}$ is contained in $\Theta\cap\Theta_x$ for
all $s$, so is
$Z_{\rm red}+A$, where $A$ is  the neutral component of the closed subgroup
generated by $x$.
It follows that $Z_{\rm red}+A=Z_{\rm red}$, hence $Z_{\rm red}$ contains a
translate $A_z$ of $A$,
which may further be assumed to satisfy $\codim _X
(A_z\cap\Sing\Theta)>3$. If $Z$ is not reduced, it is contained in the singular
locus of $\Theta\cap\Theta_x$, hence so is $A_z$.
As in the proof of proposition 2.6 of
[D],  it follows from the Jacobian criterion
and the inequality $\codim _X
(A_z\cap\Sing\Theta)>3$ that $x$ is in the kernel of the restriction
homomorphism
$\Pic ^0(X)\ra\Pic ^0(A_z)$, hence so is $A$. The composed homomorphism $A\ra
\Pic ^0(A)$
is therefore zero. Since it is the morphism associated with the restriction of
the polarization
$\lambda$ to $A$, this implies  $A=0$, which contradicts the fact that $x$ is
not torsion. Hence $Z$
is reduced.\cqfd

\ind In the case of a degenerate trisecant, this lemma allows us to prove the
following improvement
on theorem 2.2 of [D].

{\pc Theorem} 6.2.-- {\it Let $(X,\lambda)$ be a complex  \ppav , let $\Theta$
be
a symmetric representative of the polarization and let $K:X\ra
|2\Theta |^*$ be the Kummer morphism. Assume that
there exist two points $u$ and $v$ of $X$ such that:

{\parindent=1.5cm\item{\rm (i)}{\it the points $K(u)$ and $K(v)$ are distinct
and non-singular on
$K(X)$ and the line that joins them is tangent to $K(X)$ at $K(u)$,
\item {\rm (ii)} the point $u$ is not torsion,
\item {\rm (iii)}$\codim _X\bigl(\Sing\Theta\cap\bigcap_{s \in {\bf Z}}
\Theta_{2su}\bigr)>3$.}\par} Then $(X,\lambda)$ is isomorphic to the Jacobian
of a smooth
non-hyperelliptic algebraic curve.} \medskip \ind Note that  by [BD], condition
(i) implies:
$\codim _X\Sing\Theta\le 4$. \smallskip
{\bf Proof.} We keep the notation of the proof of theorem 3.1. The point is to
show that
 $R_n$ vanishes on the scheme $\Theta_u\cap\Theta_{-u}$.  Let $Z$ be a
component of
$\Theta_u\cap\Theta_{-u}$. If $Z_{\rm red}$ is {\it not\/} contained in
$\bigcap_{s \in {\bf Z}}
\Theta_{u+2su}$, lemma 3.14 implies that $R_n$ vanishes on $Z$. Otherwise,
lemma 6.1 implies that
$Z$ is reduced. On page 9 of [D], it is proved that $R_n^2$ vanishes on
$\Theta_u\cap\Theta_{-u}$.
Since $Z$ is reduced, it follows that $R_n$ vanishes on $Z$. Hence  $R_n$
vanishes on all components
of $\Theta_u\cap\Theta_{-u}$, which proves the theorem.\cqfd

\ind Lemma 6.1 has an obvious analog when $x$ is replaced by a non-zero vector
field on $X$. This
yields as above an algebraic proof of the Novikov conjecture with the sole
extra hypothesis:
$\codim _X\bigl(\Sing\Theta\cap\bigcap_{s\geq 0}D_1^s\Theta\bigr)>3$. However,
this method does not
seem to work in the non-degenerate case.

 \saut \saut
\centerline {\pc References}
\saut
\parskip=2mm
\hangindent=1cm
[AC] E. Arbarello, C. De Concini, {\it Another proof of a conjecture of S.P.
Novikov on
periods of abelian integrals on Riemann surfaces\/}, Duke Math. J. {\bf 54}
(1984),
163--178.

\hangindent=1cm
[BD] A. Beauville, O. Debarre: {\it Une relation entre deux approches du
probl\`eme de Schottky\/},
Inv. Math. {\bf 86} (1986), 195--207.

\hangindent=1cm
[D] O. Debarre: {\it Trisecant Lines And Jacobians\/}, J. Alg. Geom. {\bf 1}
(1992), 5--14.

\hangindent=1cm
[S] T. Shiota: {\it Characterization of Jacobian varieties in terms of soliton
equations
\/}, Invent. Math. {\bf 83} (1986), 333--382.

\hangindent=1cm
[W] G. Welters: {\it A criterion for Jacobi varieties\/}, Ann. of Math. (2)
{\bf 120} (1984),
497--504.

\bye